\documentclass[a4paper,11pt]{article}
\usepackage[english]{babel}
\usepackage{paralist,url,verbatim,anysize,enumitem}
\usepackage{amscd}
\usepackage{xcolor} 
\usepackage{subfig} 
\usepackage{svg} 
\usepackage{mathrsfs, mathtools}
\usepackage{cancel} 
\usepackage{xfrac,faktor}
\usepackage[e]{esvect}
\usepackage{amsmath,amsthm,amssymb,amsfonts}
\usepackage{graphicx, graphics}
\usepackage[bookmarksopen=true]{hyperref}
\hypersetup{colorlinks=true, citecolor=blue}
\usepackage{bbm}
\usepackage[font=small,labelfont=bf]{caption}
\usepackage{bm} 
\usepackage[symbol*]{footmisc} 
\usepackage{tikz} 
\usepackage[]{geometry}
\numberwithin{equation}{section} 
\numberwithin{figure}{section}

\theoremstyle{plain}
\newtheorem{theorem}{\bf Theorem}[]

\newtheorem{lemma}[theorem]{Lemma}

\newtheorem{thmnonumber}{\bf Theorem}

\theoremstyle{definition}

\newtheorem{remark}[theorem]{Remark}

\newtheorem{definition}[theorem]{Definition}

\newcommand\eqdef{\mathrel{\overset{\makebox[0pt]{\mbox{\normalfont\tiny\sffamily def}}}{=}}}

\definecolor{mypink}{RGB}{215, 5, 234}
\definecolor{lemonchiffon}{RGB}{255, 250, 205}

\renewcommand{\textdagger}{$**$}
\renewcommand{\textdaggerdbl}{$**$}
\renewcommand{\textparagraph}{$*$}
\DefineFNsymbolsTM{otherfnsymbols}{%
  \textparagraph \mathparagraph
  \textdagger    \dagger
  \textdaggerdbl \ddagger
  \textasteriskcentered *
  \textbardbl    \|%
}%
\setfnsymbol{otherfnsymbols}

\begin{document}

\author{Bruno Benedetti \\
\small Dept. of Mathematics\\
\small University of Miami\\
\small bruno@math.miami.edu
\and 
Marta Pavelka\\
\small Dept. of Mathematics\\
\small University of Miami\\
\small pavelka@math.miami.edu}

\title{\vspace{-2em}Higher-dimensional counterexamples to Hamiltonicity}
\maketitle

\begin{abstract}
For $d \ge 2$, we show that all graphs of $d$-polytopes have a Hamiltonian line graph if and only if $d \ne 3$: We exhibit a graph of a $3$-polytope on $252$ vertices whose line graph does not even have Hamiltonian paths. Adapting a construction by Gr\"unbaum and Motzkin, for large $n$ we also construct simple $3$-polytopes on $3n$ vertices in whose line graph any simple path is shorter than $10 n^{\alpha}$, for some constant $\alpha<1$.  
Moreover, we give four elementary counterexamples of plausible extensions to simplicial complexes of four famous results in Hamiltonian graph theory. 

\end{abstract}
 

\vspace{0.5 em}


\section*{Introduction}
\label{intro}

 Hamiltonian graphs are named after sir William Rowan Hamilton, who in 1857 invented a commercial puzzle of finding a cycle on the edges of a dodecahedron, with the condition that every vertex should be visited exactly once. Graphs with a path that visits all vertices once are instead called \emph{traceable}. In 1880 Tait conjectured that all graphs of simple $3$-polytopes are Hamiltonian \cite{Tait}; this was disproved 66 years later by Tutte's graph \cite{Tutte}. In 1962 Gr\"unbaum--Motzkin  \cite{GrMo} and Moon--Moser  \cite{MoonMoser} constructed $3$-polytopes whose graph is `far from being traceable'. By building a pyramid over each face in such examples in an iterative way, they obtained $d$-polytopes with non-Hamiltonian graph for every $d \ge 4$.  Barnette's 1970 conjecture \cite{Barnette3} that ``for $d \ge 4$, all graphs of simple $d$-polytopes are Hamiltonian'' is still open. In 2001 Mohar made the opposite conjecture \cite{Mohar} that for every integer $d\geq 3$, there exists a simple $d$-polytope whose graph is not Hamiltonian. Hamiltonian graphs are still at the center of interesting research, see for example \cite{LWY, PR24}.
 
The line graph of a graph $G$ has the edges of $G$ as its vertices, with adjacency determined by the adjacency of the edges in $G$. This is a classical construction in graph theory, going back to Whitney, cf.~\cite{Beineke}. In the first section of the present paper, we focus on  \emph{line-Hamiltonian} graphs, which are graphs whose line graph is Hamiltonian. All Eulerian and all Hamiltonian graphs are line-Hamiltonian \cite{Chartrand}, but the converse is false: For example, Tutte's graph mentioned above is line-Hamiltonian, even if it is neither Eulerian nor Hamiltonian. We establish the following relation to polytopality:

\begin{thmnonumber}[{Theorems \ref{prop:LineHamPolytopes}, \ref{prop:GM}}]
 Let $d \ge 2$ be an integer. All graphs of $d$-polytopes are line-Hamiltonian if and only if $d \ne 3$. In fact, there is a $3$-polytope on   $252$ vertices whose graph has a  line graph that is not traceable. Moreover, for any $n$, there exists a simple $3$-polytope on $3n$ vertices whose graph has a line graph in which every simple path is shorter than $10 n^{\alpha}$, for some constant $\alpha<1$.  
\end{thmnonumber}

Thus $3$-polytopes can be `arbitrarily far' from being line-traceable. See Theorem \ref{prop:PseudoHamStar} for an extension to spheres, pseudomanifolds, or Gorenstein subspace arrangements. En passant, we prove in Theorem \ref{prop:klee}
that the dual graph of any polytopal pseudomanifold is $(d+1)$-connected, a fact that in the simplicial case had been proven long ago by Klee \cite{Klee}.

Through extensive computer search, we also tested a big portion of the non-Hamiltonian planar 3-regular graphs generated by Brinkmann and McKay \cite{BM}. All of them turned out to be line-Hamiltonian
(Theorem \ref{thm:Computation}).

In the second part of the paper,  we turn to the goal of extending the theory of Hamiltonian graphs to simplicial complexes of higher dimensions. In 1952 Dirac \cite{Dirac} proved that any graph with $n \ge 3$ vertices is Hamiltonian if every vertex has at least $\frac {n}{2}$ neighbors. This was later strengthened by P\'osa \cite{Posa} and Chv\'atal \cite{Chvatal}, who characterized the degree sequences that force a graph to have Hamiltonian paths or cycles~\cite{Chvatal}, and managed to show that all self-complementary graphs are traceable  \cite{Chvatal}. A crucial observation in Chv\'atal's proof is that every ``maximally non-Hamiltonian'' graph admits a Hamiltonian path. 

For $d$-dimensional simplicial complexes, three distinct notions of ``tight-Hamiltonian paths'', ``loose-Hamiltonian paths'', and ``weakly--Hamiltonian paths'' were introduced in \cite{KO, KK99, BSV}. ``Tight-'' implies ``loose-'' implies ``weakly-Hamiltonian", and in general the implications are strict, except when d=1, i.e. for graphs: In this case, all three notions coincide with the usual Hamiltonicity property. 
With these notions, some Dirac-type theorems have been achieved in higher-dimensions, but only for complexes with a huge number of vertices \cite{HS10, Kand10, RSR08}, or with large ridge degree \cite{KK99}, or for complexes that are already known to be traceable~\cite{BSV}. We show a few obstacles to extending Chvátal's aforementioned observation, Dirac's theorem, or the traceability of self-complementary graphs, to higher dimensions: 

\begin{thmnonumber}[{Theorems \ref{prop:D_d}, \ref{prop:2}, \ref{prop:SC}, \ref{prop:Fl}}] We construct:
\begin{compactitem}
\item For each $d\ge 2$, a $d$-complex without (tight, loose, or weak) Hamiltonian paths that is maximally non-weakly-Hamiltonian; 
\item For each $d\ge 2$, a $d$-complex without (tight, loose, or weak) Hamiltonian paths in which any of the $n=2d+2$ vertices is in at least $ \binom{2d}{d-1}> \frac {n}{2}$ facets.
\item A self-complementary $2$-complex without (tight, loose, or weak) Hamiltonian paths.
\item A $2$-strongly connected $2$-complex whose ``square'' has no (tight, loose, or weak) Hamiltonian cycle. 
\end{compactitem}
\end{thmnonumber}

The last item is an obstacle to extending Fleischner's theorem, which claims that the ``square'' of every $2$-connected graph $G$ is Hamiltonian \cite{Flei}; see Section \ref{sec:3} for definitions.

\section{Line-Hamiltonian and non-line-Hamiltonian polytopes} 
By $\Sigma^{d}_{n}$ we mean the $d$-skeleton of the $(n-1)$-simplex on the vertices $1, \ldots, n$. By $\,H_i$ we mean the $d$-face of $\Sigma^d_{n}$ with vertices $i, i+1, i+2, \ldots, i+d$, where the sum is taken modulo~$n$. Recall that the line graph $L(G)$ of a graph $G$ is the graph with one vertex for each edge of $G$, and two vertices of $L(G)$ are joined by an edge if and only if the corresponding edges in $G$ share a vertex.

\begin{definition} A graph $G$ is \emph{line-Hamiltonian} if it satisfies any of the following equivalent conditions (see e.g. \cite{HNW} for the equivalence of (ii) and (iii)):
\label{def:cycledef}
\begin{compactenum}[\rm (i)]
\item Its edges can be labeled $1, 2, \ldots, e$, such that for all $i \in \{1, \ldots, e\}$, the edges $i$ and $i+1$ intersect; here $e+1$ is identified with $1$.
\item Its line graph is Hamiltonian.
\item Each edge has at least one endpoint on some cycle $P$ inside the graph (possibly consisting of a single vertex).
\end{compactenum}
Similarly, we call a graph \emph{line-traceable} if its line graph is traceable. A \emph{line-Hamiltonian cycle} (resp. \emph{line-Hamiltonian path}) is a cycle (resp. path) that visits every vertex of the line graph exactly once.
\end{definition}

\newpage

It is easy to see that all Eulerian and all Hamiltonian graphs are line-Hamiltonian \cite{Chartrand}. The converse is false: the graph on 6 vertices with edged $12, 23, 24, 34, 35, 45, 56$ is neither Eulerian nor Hamiltonian, yet the line graph is both. Also, all line-Hamiltonian graphs are connected, though not necessarily $2$-(edge)-connected, as shown by any cycle with some pendants attached. Note that if we attach a two-edge path to a cycle, we obtain a non-line-Hamiltonian graph. 
 
More interesting examples arise from polytope theory. A $d$-polytope is the convex hull of finitely many points in $\mathbb R^d$ that do not belong to a common hyperplane. A $d$-polytope is called simple if every vertex is contained in exactly $d$ facets. In 1946 Tutte constructed the first example of a simple 3-polytope whose graph is not Hamiltonian \cite{Tutte}, thereby disproving a long-standing conjecture by Tait \cite{Tait}. It turns out that Tutte's graph, as well as many others alike, is line-Hamiltonian:
 
  \begin{theorem} \label{prop:TutteLH}
 The following non-Hamiltonian graphs of simple $3$-polytopes are line-Hamiltonian: Tutte's graph, the Lederberg--Bosák--Barnette graph, the Faulkner--Younger graphs on $42$ and $44$ vertices, the Grinberg graphs on $42$, $44$, and $46$ vertices, and the Thomassen graph on $94$ vertices.
\end{theorem}

 
 \vskip-3mm
 \begin{figure}[htb] 
    \centering
            \subfloat[]{\includegraphics[width=9.5em]{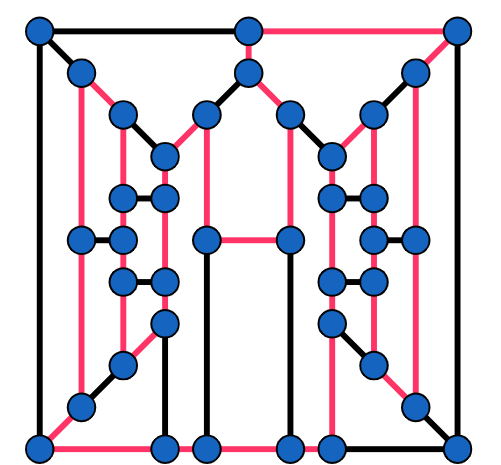}
		\label{fig:Barnette}}
   		\hfill
    	\subfloat[]{\includegraphics[width=9em]{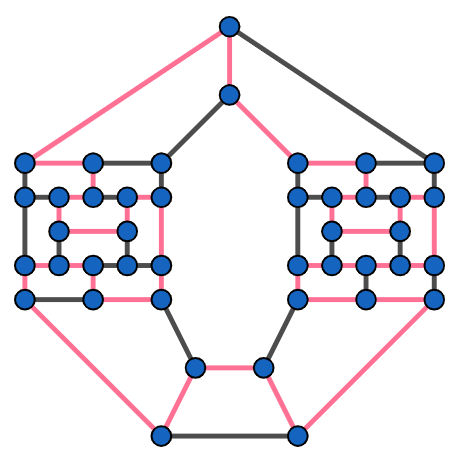}
		\label{fig:Faulkner}}
    	\subfloat[]{\includegraphics[width=9em]{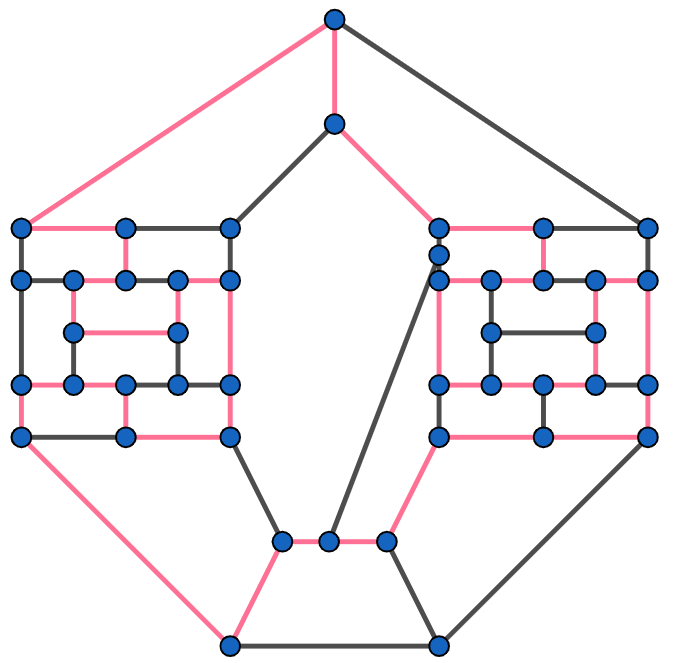}
		\label{fig:Faulkner44}}
  		\hfill
    	\subfloat[]{\includegraphics[width=9em]{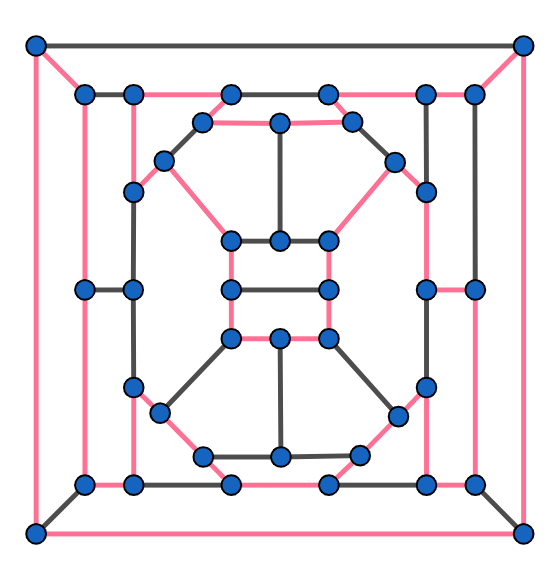}
		\label{fig:Grinberg42}}
            \hfill
    	\subfloat[]{\includegraphics[width=9em]{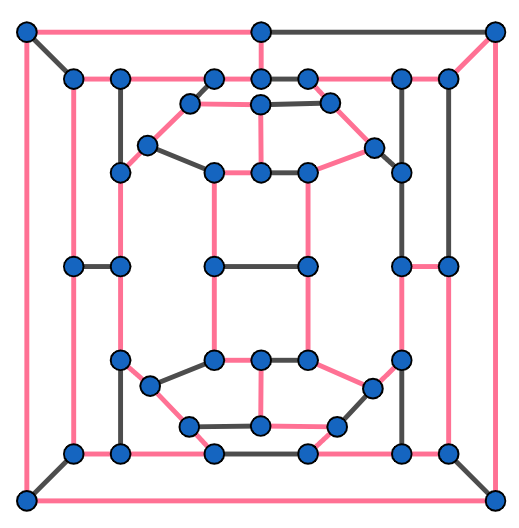}
		\label{fig:Grinbeg44}}
  		\hfill
    	\subfloat[]{\includegraphics[width=15em]{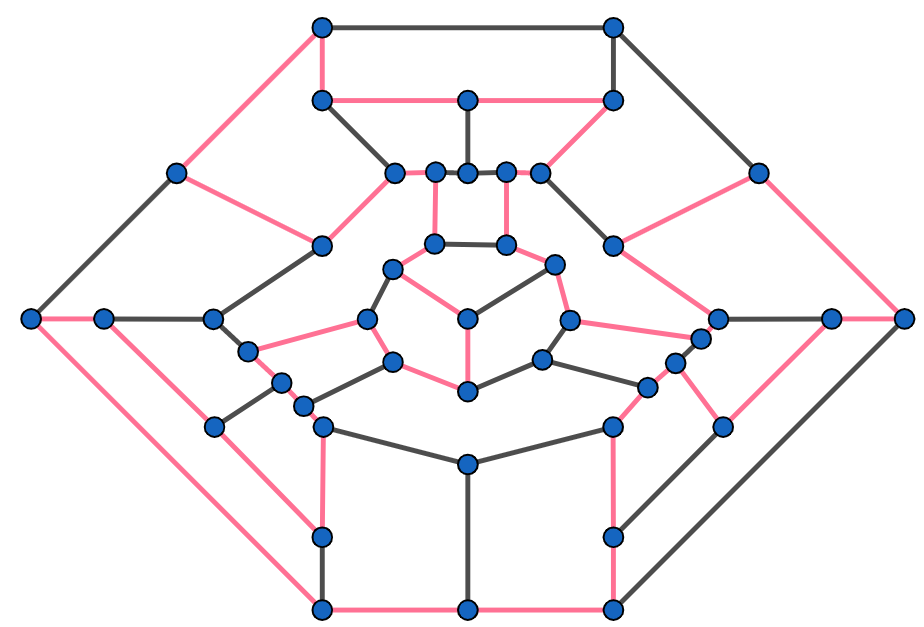}
		\label{fig:Grinbeg46}}
            \hfill
    	\subfloat[]{\includegraphics[width=9.5em]{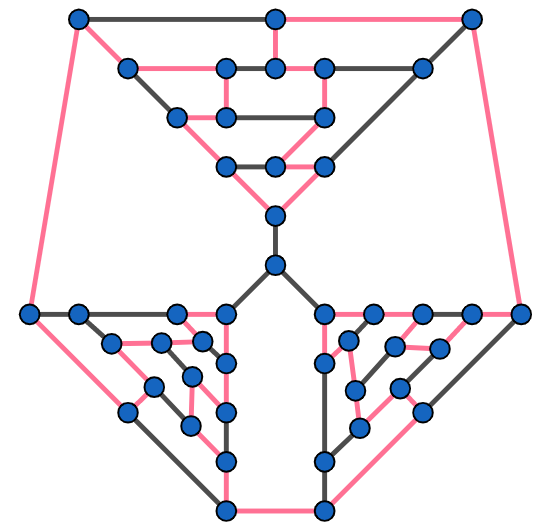}
		\label{fig:Tutte}}
  		\hfill
      	\subfloat[]{\includegraphics[width=14em]{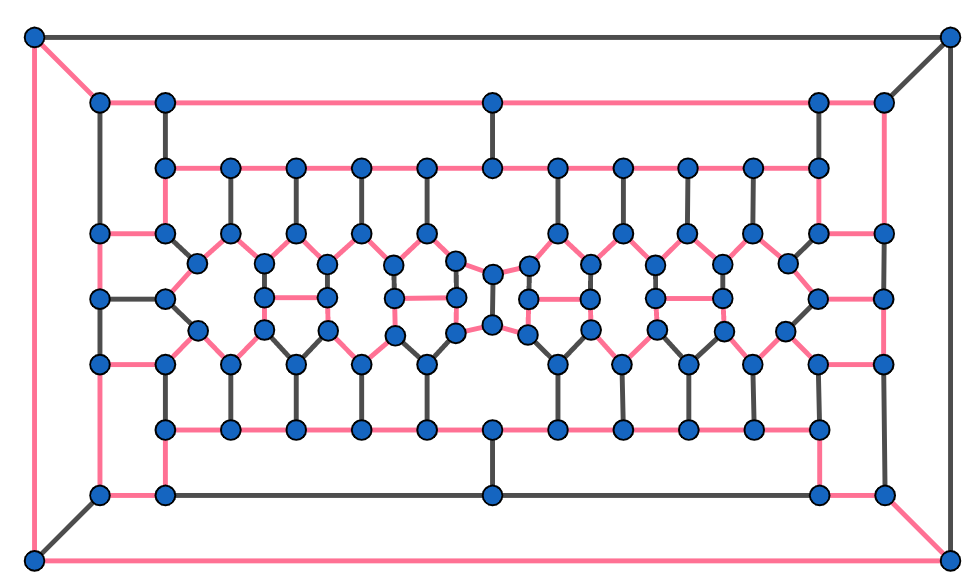}
		\label{fig:Thomassen}}
            \caption{Cycles proving the line-Hamiltonian property for: (a) the Lederberg--Bosák--Barnette graph; (b),(c) the Faulkner--Younger graph on 42 resp.~44 vertices; (d), (e), (f) the Grinberg graph on 42, resp.~44, resp.~46 vertices; (g) Tutte's graph; (h) the Thomassen graph on 94 vertices.}
		\label{fig:examples}
\end{figure}

\begin{proof}[Proof by picture]
For each of the listed non-Hamiltonian regular graphs, in Figure \ref{fig:examples} we exhibit a cycle (in red) that via Definition \ref{def:cycledef}, part (iii), proves the line-Hamiltonian property. 
\end{proof}

\newpage

\begin{theorem} \label{thm:Computation}
All of the following plane $3$-regular graphs are line-Hamiltonian:  
\begin{compactitem}
\item[--] those with no faces of size $3$ and with up to $48$ vertices;
\item[--] those that are cyclically $4$-connected with up to $50$ vertices;
\item[--] those with no faces of size $3$ or $4$, with cyclic connectivity of exactly $4$, and with up to $70$ vertices;
\item[--] those that are cyclically $5$-connected with up to $74$ vertices.
\end{compactitem}
\end{theorem}

\begin{proof}
   By an extensive computer search \cite{Pav23} on the non-Hamiltonian graphs generated by Brinkmann and McKay \cite{BM}, we verified that all the graphs listed above are line-Hamiltonian. The class with no faces of size 3 or 4 and with cyclic connectivity of exactly 4 was the most challenging one for our server. Among these graphs with 70 vertices, roughly a dozen of them each required two weeks to be verified. Our script is run on a Quad Core CPU with the Neoverse N1 model, operating at a speed of 3.0 GHz.
\end{proof}

This suggests an update to Tait’s question: are perhaps all polytopal graphs \emph{line}-Hamiltonian? The answer turns out to be negative for $d=3$, and positive otherwise:

\begin{theorem} 
\label{prop:LineHamPolytopes}
Let $d \ge 2$ be an integer. The graphs of all $d$-polytopes are line-Hamiltonian if and only if {$d \ne 3$}. In fact, there is a simple $3$-polytope with $124$ vertices whose graph is  line-traceable but not line-Hamiltonian, and a simple $3$-polytope on $252$ vertices whose graph is not line-traceable.
\end{theorem}

\begin{figure}[t] 
\vspace{-2em}
    \centering
            \subfloat[]{\includegraphics[width=12em]{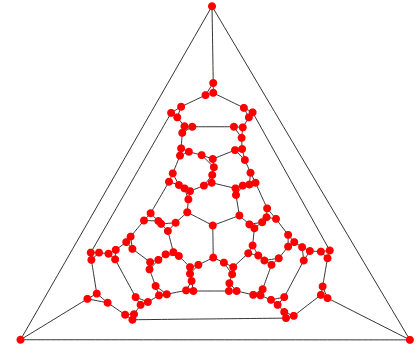}
		\label{fig:Grum1}}
  		\hfill
    	\subfloat[]{\includegraphics[width=12em]{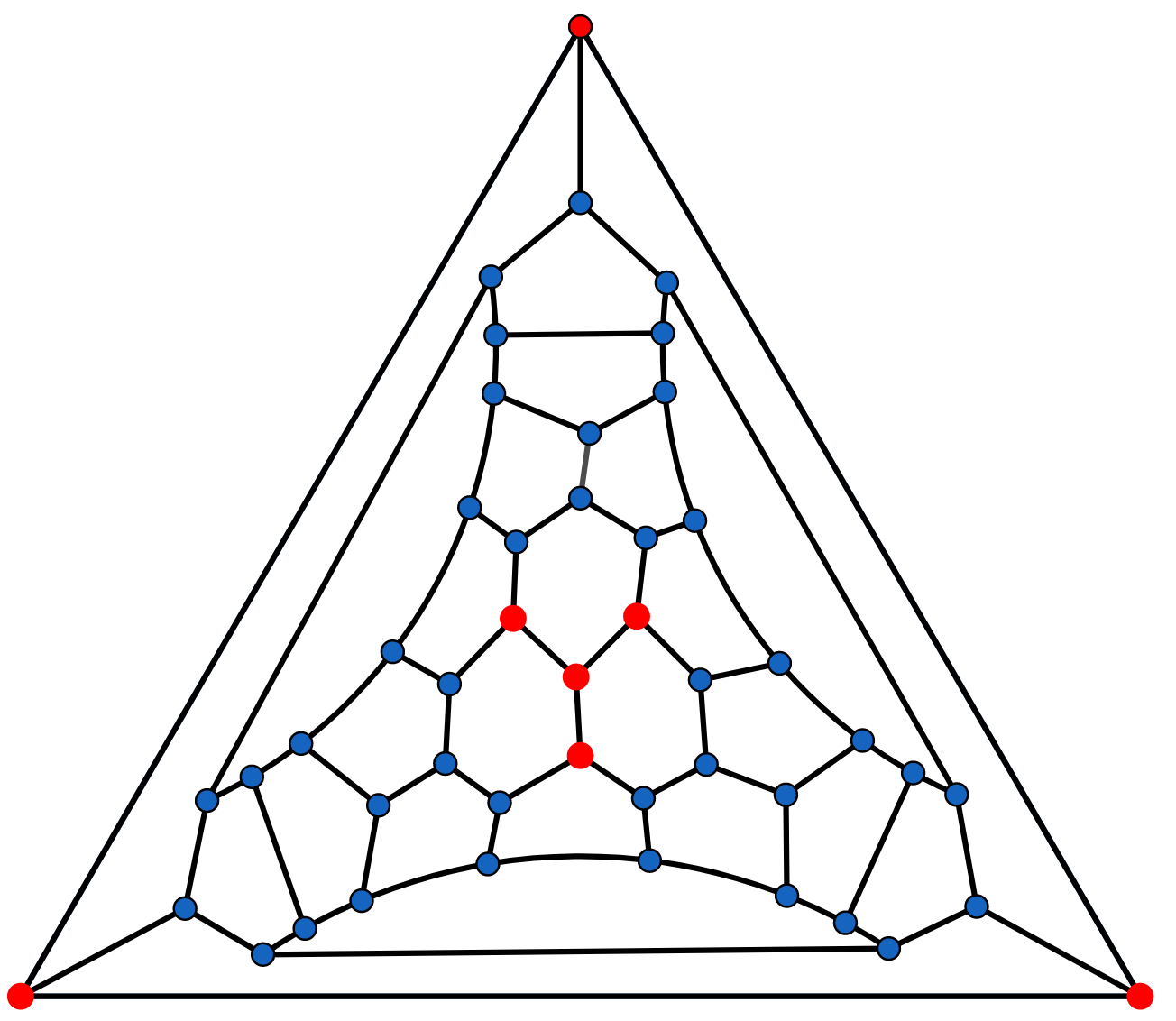}
		\label{fig:Grum2}}
            \hfill
    	\subfloat[]{\includegraphics[width=12.5em]{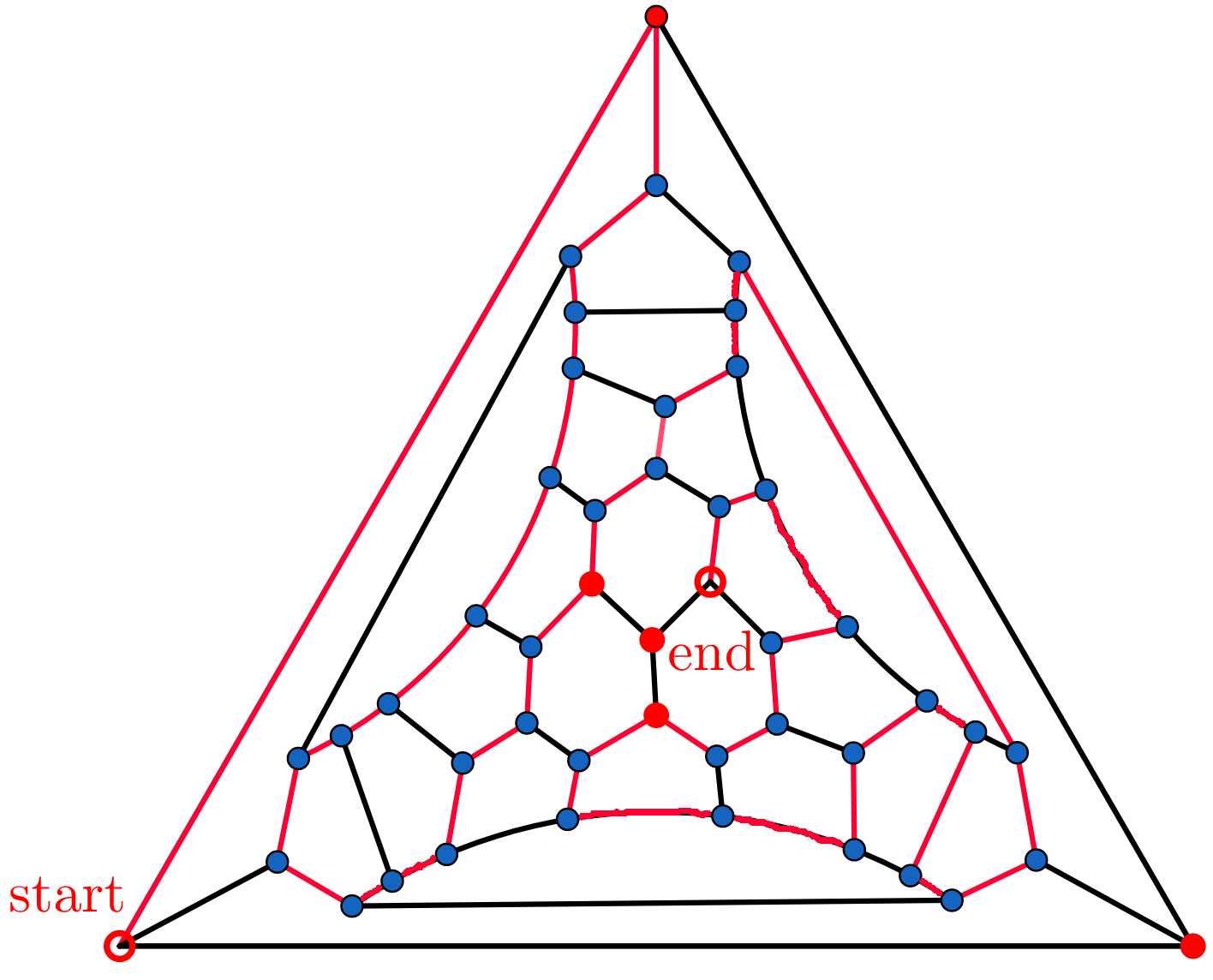}
		\label{fig:GrumLineTrace}}
   		\caption{(a) The planar 3-regular Gr\"{u}nbaum graph of a simple 3-polytope that is not Hamiltonian, nor line-Hamiltonian. (b) Simplification of the Gr\"{u}nbaum graph. (c) Spanning path of the Gr\"{u}nbaum graph implying existence of a line-Hamiltonian path.}
		\label{fig:Grumm}
\end{figure}

\begin{figure}[t] 
\vspace{-4em}
    \centering
    	\subfloat[]{\includegraphics[width=11em]{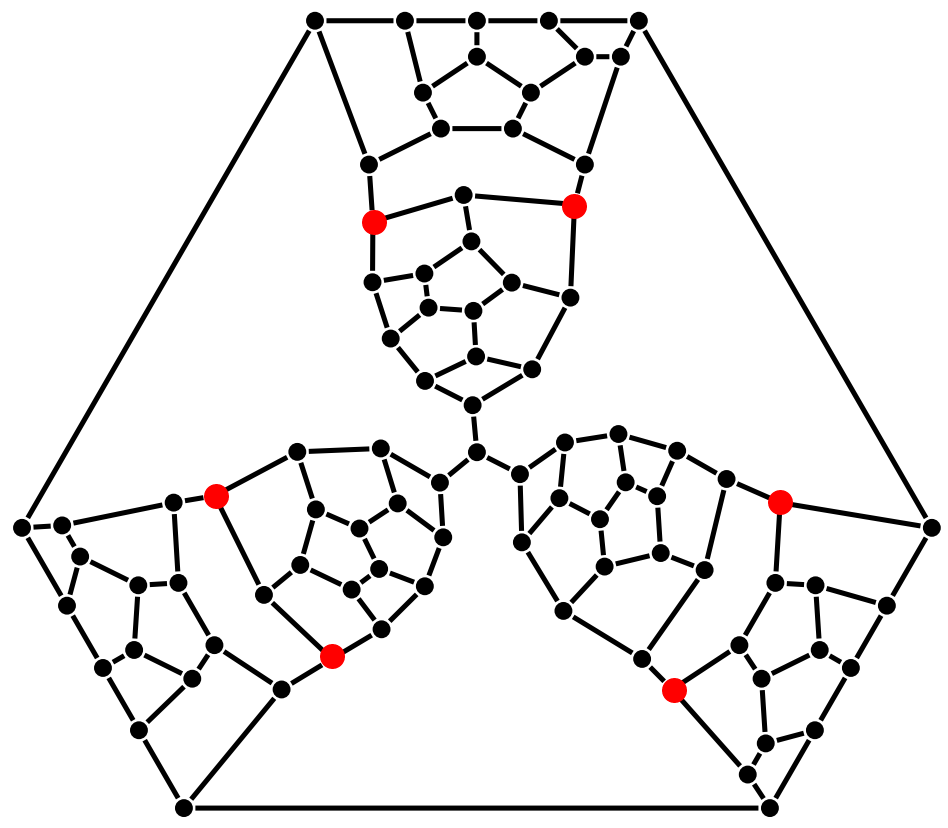}
		\label{fig:NonTraceable}}
            \hfill
    	\subfloat[]{\includegraphics[width=10em]{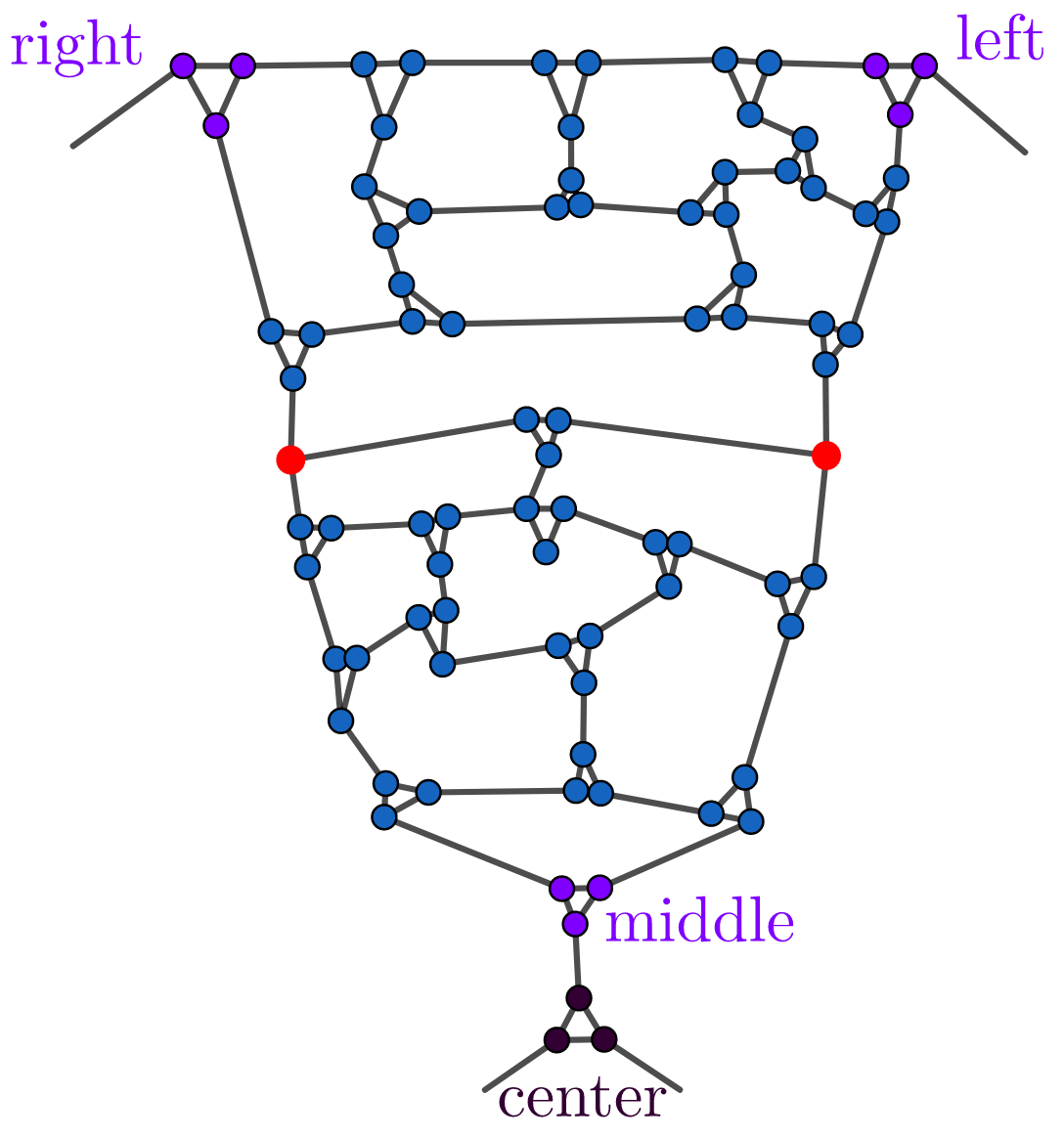}
		\label{fig:NonTraceableBlock}}
		\caption{(a) Zamfirescu's $3$-regular graph is not traceable, but it is line-traceable. It can be modified by replacing all but the six red vertices with triangles into a non-line-traceable graph of a simple 3-polytope with 252 vertices. (b) A fragment of the simple 3-polytope with non-line-traceable graph with 252 vertices.}
		\label{fig:Grum}
\end{figure}

\begin{proof} By Balinski's theorem, graphs of $d$-polytopes are $d$-connected, and so are their line graphs. Chen, Lai, Lai, and Weng \cite{CLLW} proved in 1993 that if a connected graph $G$ does not have vertices of degree 3 and its line graph $L(G)$ is $4$-connected, then $L(G)$ is Hamiltonian.  This immediately implies that the graphs of all $d$-polytopes are line-Hamiltonian for $d \ge 4$. (The same result would also follow from a 2012 paper of Broersma, Ryjáček, and Vrána \cite[Corollary 30]{BRV}, and for polytopes of dimension $5$ or more, from Kaiser-Vrána's 2012 result that every $5$-connected line graph with minimum degree $\ge 6$ is Hamiltonian \cite{KV}, as the line graph of the graph of a $d$-polytope has minimum degree $\ge 2d$.) 

The smallest polytopal graph for which we were able to exclude the line-Hamiltonian property, is the planar 3-regular Gr\"{u}nbaum graph $H$ in Figure \ref{fig:Grum1}, on 124 vertices.  Every line-Hamiltonian cycle has to visit any of the small triangles exactly once since the graph is 3-regular. The graph $G$ in Figure \ref{fig:Grum2} is constructed from $H$ by replacing each of the small triangles with a blue vertex; the remaining vertices stay and are colored red. Thus $G$ is line-Hamiltonian if and only if there is a Hamiltonian subgraph in $H$ that contains all the blue points and a suitable combination of the red points. However, with a computer check via \textsc{SageMath} \cite{Sage}, we verified that no such subgraph of $H$ exists. Thus $G$ is not line-Hamiltonian; Figure \ref{fig:GrumLineTrace} shows however that it is line-traceable.

It remains to construct a simple 3-polytope that is neither line-Hamiltonian nor line-traceable. In the hope of getting an example as small as possible, we start with the smallest known non-traceable graph  of a simple 3-polytope; that is the Zamfirescu's graph $Z$ on 88 vertices from Figure \ref{fig:NonTraceable}. We replace all but the six red vertices with triangles and denote this new graph $\tilde Z$. We claim that this cannot be line-traceable. In fact, each of the three fragments has only three edges coming out of it: toward the left, right, and middle; as illustrated in Figure \ref{fig:NonTraceableBlock}. Hence each fragment can be passed through only once by any spanning path of $\tilde Z$. Moreover, the center triangle of $\tilde Z$ can be passed through only once. From these observations, it is clear that any spanning path of $\tilde Z$ has to pass through at least one of the fragments in one go, either from left to right or from the middle to right. A spanning path of $\tilde Z$ would then correspond to a simple path in $Z$ that meets every black vertex exactly once and may or may not pass through the red vertices. However, with a computer check via \textsc{SageMath} \cite{Sage}, we verified that there is no spanning path of the fragment going from left to right nor from the middle to right, regardless of which of the red vertices we include.
\end{proof}

In the early Sixties, Grünbaum--Motzkin \cite{GrMo} and Moon--Moser \cite{MoonMoser} constructed two families of simple $3$-polytopes whose graphs are ``arbitrarily far’’ from being traceable, in the sense that the longest path in these graphs involves only a fraction $f(n)$ of the $n$ vertices that tends to zero for $n$ large. Our next goal is to strengthen this classical result by replacing “traceable” with “line-traceable” in the conclusion. This is not obvious, as the next theorem shows:

\begin{theorem}
    The graphs of all polytopes from the Moon--Moser family \cite{MoonMoser} are line-Hamiltonian. 
\end{theorem}

\vspace{-1em}

\begin{figure}[h] 
    \centering
            \subfloat[]{\includegraphics[height=7.5em]{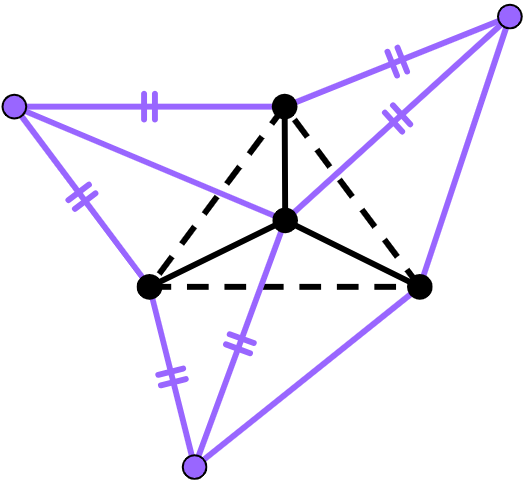}
		\label{fig:MM1}}
  		\hfill
    	\subfloat[]{\includegraphics[height=7.5em]{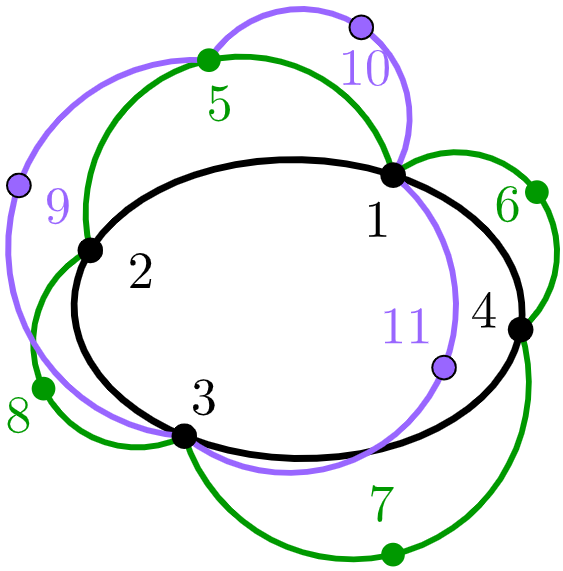}
		\label{fig:MM2}}
            \hfill
    	\subfloat[]{\includegraphics[height=7.5em]{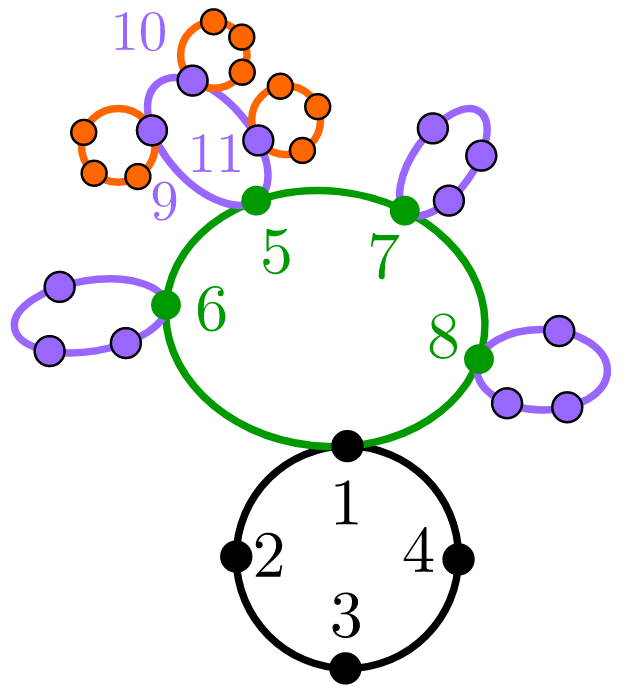}
		\label{fig:MM3}}
   		\caption{The Moon--Moser construction: (a) Zooming in on one of the $(k-1)$-st stage tetrahedra, whose vertices and edges are in black. The $k$-th stage vertices and edges are in purple. A cycle meeting all the new vertices that uses only new edges and passes through the vertex common to all three new tetrahedra is indicated by dashes on its edges. (b) Inductively, we get cycles at each stage. 
     Here we show cycles in stages 0, 1, and one cycle (out of four) 
 in stage 2. (c) To illustrate the spanning cycle, we remove all unnecessary intersections from the picture.}
		\label{fig:MoonMoser}
\end{figure}

\begin{proof}
Let us briefly recall the Moon--Moser construction. We start with a tetrahedron $G_0$, and call its vertices and edges ``the 0th stage vertices and edges''. To each boundary triangle of $G_0$, we attach a new tetrahedron by one of its boundary triangles. Let $G_1$ be  the resulting polytope $G_1$; it has 12 boundary triangles, 4 new vertices (which we call ``1st stage vertices''), and 12 new edges (the ``1st stage edges''). We repeat the procedure until we reach, say, $G_k$. In the graph of $G_k$, as we zoom in into one of the tetrahedra from the stage $(k-1)$, say $\Delta$, we see the three new tetrahedra attached to it, see Figure \ref{fig:MM1}. In this subgraph, we can always find a cycle visiting all its vertices of the $k$-th stage using only the edges of the $k$-th stage and passing through the vertex of $\Delta$ where all the three attached tetrahedra meet. Since a spanning cycle can pass through a vertex multiple times (depending on its degree), from these small cycles we get inductively a spanning cycle of the whole graph, see Figures \ref{fig:MM2}, \ref{fig:MM3}.  This shows that all Moon--Moser graphs are line-Hamiltonian, even if they are far from being traceable.
\end{proof}

Nonetheless, we will now show that a suitable modification of the Grünbaum-Motzkin construction succeeds in producing simple polytopes whose graphs are arbitrarily far from being line-traceable.

\begin{theorem} \label{prop:GM} 
    For any positive integer $n$, there is a graph $G^{\Delta}$ of a simple $3$-polytope on $3n$ vertices with $p(L(G^{\Delta}))< 10 n^{\alpha}$, where $\alpha<1$ (e.g., $\alpha = 1-2^{-19}$).
\end{theorem}

\begin{figure}[h] 
    \centering
            \subfloat[]{\includegraphics[height=6em]{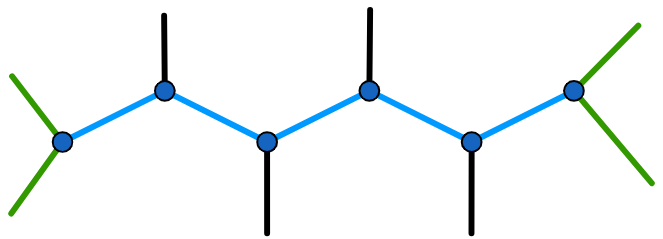}
		\label{fig:GrunPath1}}
  		\hfill
    	\subfloat[]{\includegraphics[height=7em]{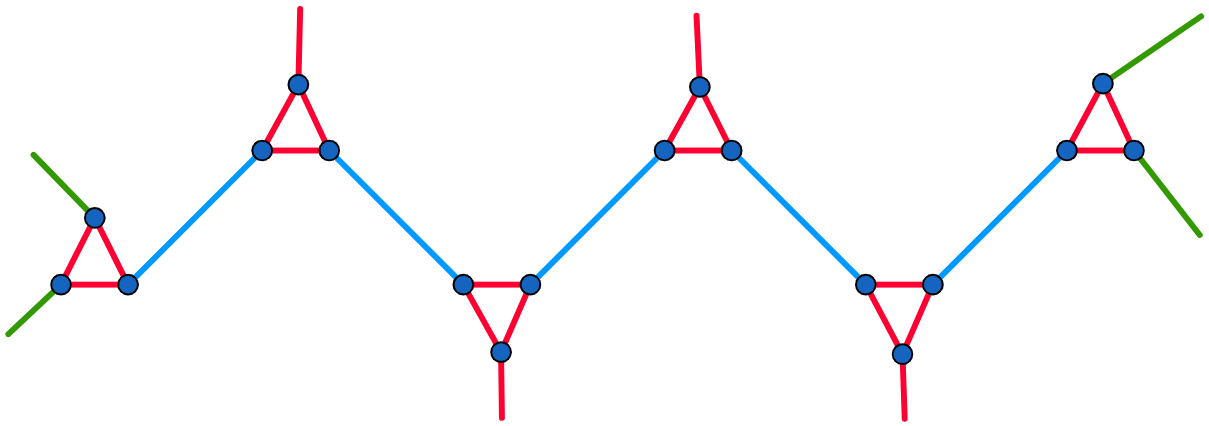}
		\label{fig:GruPath2}}
   		\caption{ (a) Blue edges represent a path $p$ of maximal length in a Gr\"unbaum--Motzkin graph $G_n$. The green edges cannot lead to a new vertex, or else $p$ would not be maximal. (b) A corresponding path in $G^{\Delta}_n$ spans a subgraph with $\le 4$ edges for each interior vertex of $p$ (in red), plus the original $m-1$ edges (in blue), plus 2 edges for the boundary vertices (in red).}
		\label{fig:GrunMotz}
\end{figure}

\begin{proof}
    Following the notation from \cite{GrMo}, let $p(G)$ denote the maximal length of a simple path in $G$. For fixed $n$, let us denote by $G_n$ the Gr\"unbaum--Motzkin graph of a simple 3-polytope with $n$ vertices for which $p(G_n)<2n^{\alpha}$ from \cite{GrMo}. Since each vertex has degree 3, we can replace each vertex of $G_n$ with a triangle without destroying the 3-regularity. Let us denote the resulting graph by $G^{\Delta}_n$. Since $G^{\Delta}_n$ is 3-regular, every path (even not a simple one) can pass through each triangle only once. Now, say we have a path $p$ in $G_n$ of length $p(G_n)=m$. Any corresponding path in $G^{\Delta}_n$ spans a subgraph with at most 4 edges for each interior vertex of $p$, the original $(m-1)$ edges, and 2 edges for the boundary vertices (edges of the boundary vertices cannot lead to a new vertex, else we could extend $p$), see Figure \ref{fig:GrunMotz}. Therefore,
    \begin{equation*}
        p(L(G^{\Delta}_n)) \leq 4(m-2) + (m-1) + 2 \cdot 3 =5m-3<5(2n^{\alpha}) -3 < 10n^{\alpha} \ .  
    \end{equation*}
\end{proof}

As a conclusion to this section, we should mention that both Balinski's theorem and its equivalent dual formulation (``\emph{dual} graphs of $d$-polytopes are $d$-connected'') have been extended beyond the world of polytopes:

\begin{definition} A \emph{$d$-pseudomanifold} is a strongly-connected $d$-dimensional polytopal complex where every $(d-1)$-face lies in exactly $2$ facets.
\end{definition} 

\begin{definition} Let $P$ be a $d$-pseudomanifold. By a \emph{path of $d$-faces connecting $F'$ and $F''$ in $P$} we mean a sequence $F_1, \ldots, F_s$ of $d$-faces of $P$ such that $F_1=F'$, $F_s=F''$, and for each $i \in \{1, \ldots, s-1\}$ the intersection $F_i \cap F_{i+1}$ is $(d-1)$-dimensional. We say that two paths of faces connecting $F'$ and $F''$ are \emph{interior-disjoint} if they do not have $d$-faces in common other than $F'$ and $F''$. Finally, a \emph{cycle of $d$-faces} is a path where first and last facet coincide.
\end{definition}
 
\begin{lemma}[{Barnette, \cite[Lemma 1]{Barnette}}] \label{lemma9} Let $P$ be a $d$-pseudomanifold. Let $S$ and $T$ be two $d$-faces that share a $(d-2)$-face $\sigma$. Then in $P$ there is a cycle of $d$-faces all containing $\sigma$. 
\end{lemma} 

\begin{theorem} \label{prop:klee} The graph and the dual graph of every $d$-pseudomanifold  are both $(d+1)$-connected. 
\end{theorem}

Several versions of Theorem \ref{prop:klee} were proven in the Seventies. That graphs of $d$-pseudomanifolds are $(d+1)$-connected was  first proven in 1973 by Barnette \cite{Barnette}. In 1975 Klee proved the dual result, but only for \emph{simplicial} $d$-pseudomanifolds \cite{Klee};  Klee's argument does not extend beyond the simplicial case. In 1983 Barnette proved the dual result for  homology manifolds, simplicial or not \cite{Barnette1}, but his argument via manifold duality does not extend to all pseudomanifolds. Since we were unable to find in the literature a proof of Theorem \ref{prop:klee} valid in full generality, for completeness, we provide one below. Given a $d$-dimensional simplicial complex $G$ in which every $(d-1)$-face belongs to at most two $d$-faces, we define the ``boundary'' $\partial G$ of $G$ to be the smallest subcomplex containing all $(d-1)$-faces that belong to exactly one $d$-face of~$G$.

\begin{proof}[\textbf{Proof of Theorem \ref{prop:klee}}]
Fix a $d$-pseudomanifold. Its graph is $(d+1)$-connected by the work of Barnette \cite{Barnette}, so let us focus on its dual graph. 
Let $X$ be a collection of at most $d$ facets in the $d$-pseudomanifold. Let $F'$ and $F''$ be two facets not in $X$. Let $F'=F_1,   F_2,  \ldots  ,   F_s=F''$
be a path of $d$-faces that contains as few elements of $X$ as possible. We claim that such a path must avoid $X$ altogether. In fact, suppose not. Let $G$ be the first $d$-face from $X$ that we encounter along the path. Note that the path features three consecutive $d$-faces $F, G, H$, with $F \notin X$ and $G \in X$. 
Since the dual graph of $G$ is $d$-connected (by Balinski’s theorem), in $\partial G$ between $\sigma \eqdef (F \cap G)$ and $\tau \eqdef (G \cap H)$ there are $d$ interior-disjoint paths of $(d-1)$-faces. But any path of $(d-1)$-faces in $\partial G$ is the restriction to  $\partial G$ of a sequence of facets of $P$ adjacent to $G$, each intersecting the previous one in codimension $1$ or $2$. 
We claim that any such sequence can be completed to a path of $d$-faces connecting $F$ and $H$. In fact, let $S$ and $T$ be any two consecutive $d$-faces in the sequence with $\dim (S \cap T)= d-2$. By Lemma \ref{lemma9} there is a cycle of $d$-faces that all contain $S \cap T$; adding to our sequence the portion of the cycle that does not contain $G$, we can `fill in the gap' between $S$ and $T$ in the dual graph of $P$. So the claim is proven. Moreover, all $d$-faces we insert when we fill in gaps have $(d-2)$-dimensional intersection with $G$. In particular, no face $U$ can be used to fill in two different gaps between different pairs of facets, or else $U \cap G$ would contain two distinct $(d-2)$-faces, and thus it would have to be $(d-1)$-dimensional. We conclude that there exist $d$ interior-disjoint paths of $d$-faces connecting $F$ to $H$ and avoiding $G$. But  since $|X \setminus \{G\}| < d$, at least one of these $d$ paths does not contain elements of $X$ other than (possibly) $H$. Using this path to bypass $G$, we can modify $F_1,   F_2, \ldots ,   F_s$ to obtain another path of $d$-faces that features fewer elements of $X$. This contradicts how the path $F_1,   F_2, \ldots ,   F_s$ was chosen. 
\end{proof}

\begin{theorem}
\label{prop:PseudoHamStar} Let $d$ be a positive integer. \begin{compactitem}
\item All graphs and all dual graphs of $d$-pseudomanifolds are line-Hamiltonian if and only if $d \ne 2$. 
\item All graphs of flag simplicial $d$-pseudo\-manifolds are line-Hamiltonian. 
\item All dual graphs of subspace arrangements defined by an ideal $I$ such that  $S/I$ is Gorenstein and has regularity $d+1$ are line-Hamiltonian if $d \ge 3$, and the bound is sharp.
\end{compactitem}
\end{theorem}

\begin{proof} Every $1$-pseudomanifold is a cycle, so both its graph and its dual graph are trivially Hamiltonian. By Theorem \ref{prop:klee}, graphs and dual graphs of $d$-pseudomanifolds are $(d+1)$-connected. Athanasiadis \cite{Athanasiadis} proved that the graph of every flag simplicial $d$-pseudomanifold is $2d$-connected, for $d \ge 2$; cf.~also \cite{BjVo}. Benedetti--Varbaro \cite{BV} showed that if  $I \subseteq S=\mathbb{K}[x_1, \ldots, x_n]$ is an ideal defining a subspace arrangement, with $S/I$ Gorenstein of Castelnuovo--Mumford regularity $d+1$, the dual graph of the arrangement is $(d+1)$-connected. These four results, paired with Chen--Lai--Lai--Weng's theorem  \cite{CLLW}, imply the line-Hamiltonian property as in the proof of Theorem \ref{prop:LineHamPolytopes}, if $d$ is within the bounds described for each case.  As for the sharpness of such bounds: the Gr\"{u}nbaum graph $H$ from Figure \ref{fig:Grum1} is the graph of a simple $3$-polytope $P$, whose boundary complex is a polyhedral 2-sphere with graph $H$; the boundary complex $C$ of the dual polytope of $P$ is a simplicial $2$-pseudomanifold with dual graph $H$; if $I_C$ is the Stanley--Reisner ideal of $C$, then $S/I_C$ is Gorenstein of regularity $3$, and the dual graph to $I_C$ is still $H$ \cite{BV}.
\end{proof}

\section{Four examples of non-Hamiltonian complexes} \label{sec:3}
There are three ways in the literature to extend to $d$-dimensional complexes the notion of ``Hamiltonian graph''; we recall them below:

\begin{definition}
Let $\Sigma^{d}_{n}$ be the $d$-skeleton of the $(n-1)$-simplex, with vertices labeled by $1, \ldots, n$. Let $\,H_i$ be the $d$-face of $\Sigma^d_{n}$ with vertices $i, i+1, i+2, \ldots, i+d$, where the sum is taken modulo~$n$. We say that a $d$-dimensional simplicial complex $\Delta$ on vertex set $\{1, \ldots, n\}$
\begin{compactitem}
\item \emph{admits a tight-Hamiltonian path} (resp. \emph{cycle}) if it has a vertex labeling such that all of $H_1, \dots, H_{n-d}$ (resp. all of $H_1, \dots, H_{n}$) are in $\Delta$. 
\item \emph{admits a weakly-Hamiltonian path} (resp. \emph{cycle}) if it has a vertex labeling such that $\Delta$ contains faces $H_{i_1}, \dots, H_{i_k}$ from the set $\{H_1, \dots, H_{n-d}\}$ that altogether cover all the vertices, and such that $H_{i_j}$ is incident to $H_{i_{j+1}}$ for each $j\in \{1, \dots, k-1\}$ (respectively: for each $j\in \{1, \dots, k\}$, where by $i_{k+1}$ we mean $i_1$). 
\item \emph{admits a loose-Hamiltonian path} (resp. \emph{cycle}) if it admits a weakly-Hamiltonian path (resp. cycle) with all of the intersections $H_{i_j} \cap H_{i_{j+1}}$ consisting of a single point, except possibly for one exception, where the intersection may be larger. 
\end{compactitem}
We say that $\Delta$ is \emph{maximally non-tight-Hamiltonian} (resp. \emph{maximally non-weakly-Hamiltonian}) if it does not admit any tight-Hamiltonian (resp. weakly-Hamiltonian) cycles, yet for any ``missing $d$-face'' $\sigma \in \Sigma^{d}_n$, $\sigma \notin \Delta$, the $d$-complex $\Delta' = \Delta \cup \sigma$ does admit some tight-Hamiltonian (resp. weakly-Hamiltonian) cycle.
\end{definition}

The last definition echoes the following classical one: a graph $G$ is \emph{maximally non-Hamiltonian} if $G$ is not Hamiltonian, yet for any two vertices $x,y$ not adjacent in $G$, the graph $G$ + $[x,y]$ is Hamiltonian. Any maximally non-Hamiltonian graph has a Hamiltonian path. This trick is used in the proofs of Pósa's \cite{Posa} and Chv\'atal's theorem \cite{Chvatal}; see also  \cite{CE1, CE2, Kronk, Roldugin} for more examples of that technique. The next result shows that this trick no longer works in higher dimensions:

\begin{theorem}  \label{prop:D_d} For any integer $d \ge 2$, there is  a $d$-dimensional simplicial complex $D_d$ without (tight, loose, or weak) Hamiltonian paths in which any of the $n=2d+2$ vertices is in at least $ \binom{2d}{d-1}> \frac {n}{2}$ facets. Moreover, $D_2$ is maximally non-weakly-Hamiltonian. 
\end{theorem}

\begin{proof} Let $D_d = v \ast \Sigma^{d-1}_{2d+1}$. 
This $D_d$ is a $d$-complex on $n = 2d+2$ vertices. The vertex $v$ is present in all facets, by construction. Any other vertex is in $\binom{2d}{d-1}$ facets, which is larger than $d+1$ for $d \ge 2$. So every vertex has degree $> \frac {n}{2}$. Suppose by contradiction that $D_d$  admits a weakly-Hamiltonian path or cycle. Such path/cycle must be formed by facets $H_i$'s that all contain the vertex $v$. At the same time, if $k \in \{1, \ldots, n\}$ is the label assigned to vertex $v$, the only facets $H_i$ that contain $k$ are $H_k, H_{k-1}, \ldots, H_{k-d}$: These are $d+1$ faces that altogether cover $2d+1$ vertices, namely, those with a label between $k-d$ and $k+d$. But our complex $D_d$ has $n=2d+2$ vertices, so at least one of its vertices is not covered by the weakly-Hamiltonian path/cycle. A contradiction. Finally, we verified with the software \cite{Pav21} that adding any `missing triangle' to $D_2$ produces a weakly-Hamiltonian cycle. 
\end{proof}

The computational approach above, which proved $D_d$ maximally non-weakly-Hamiltonian for $d=2$, is of course not available for larger or arbitrary $d$. The next theorem however showcases an even simpler example that works in all dimensions:

\begin{theorem} \label{prop:2}  Let $d$ be a positive integer. Any $d$-complex that is maximally non-tight Hamiltonian must admit tight Hamiltonian paths; if $d \ge 2$, it must also admit loose-Hamiltonian cycles. In contrast, for any integer $d \ge 2$, there is a $d$-dimensional simplicial complex without (tight, loose, or weak) Hamiltonian paths that is maximally non-weakly-Hamiltonian.
\end{theorem}

 \begin{proof} For the first part, let $C$ be a $d$-dimensional simplicial complex on $n$ vertices that is maximally non-tight-Hamiltonian. Let  $\Delta$ be any face of $\Sigma^d_n$ that is not in $C$. Then $C \cup \Delta$ contains some tight-Hamiltonian cycle, in which $\Delta$ necessarily appears. Let $\Delta, F_1, \ldots, F_m, \Delta$ be any such cycle. Then $F_1, \ldots, F_m$ is a tight-Hamiltonian path in $C$. For the second part, let $W^d$ be the join of $3$ disjoint points with the $(d-1)$-simplex. It is easy to see that the $d$-complex $W^d$ does not have (tight, loose, or weak) Hamiltonian paths, cf.~\cite{BSV}. Let us prove that $W^d$ is maximally non-weakly-Hamiltonian. We label the vertices of the $(d-1)$-simplex by $1, \ldots, d$. Let us call ``apices'' the vertices $d+1, d+2$, $d+3$. Let $F \in \Sigma^d_{d+3}$. There are 3 cases:
\begin{compactenum}[ (i)]
\item If $F$ contains exactly $1$ apex, then $F$ already belongs to $W^d$. 
\item If $F$ contains exactly $2$ apices, then up to permuting the labels of the first $d$ vertices and the labels of the last $3$ we can assume that
$F=[1, 2, \ldots, d-1, d+2, d+3]$.
Hence $W^d \cup F$ admits a weak Hamiltonian cycle formed by the three faces $H_1, H_{d+3}, F$.

\item If $F$ contains $3$ apices, up to relabeling 
$F=[1, 2, \ldots, d-2, d+1, d+2, d+3]$. (When $d=2$, $F=[d+1, d+2, d+3]$.) So $W^d \cup F$ admits the weak cycle
$H_1, H_{d+3}, F$.

\end{compactenum}
\end{proof}

Our next example shows that Chv\'atal's famous result that ``all self-complementary graphs are traceable'' \cite{Chvatal} does not extend to higher dimensions either.

Recall that a pure $d$-dimensional complex $S$ on $n$ vertices is called \emph{self-complementary} if it is combinatorially equivalent to the pure $d$-complex whose facets are the $d$-faces of $\Sigma^d_{n}$ that are not in $S$. 

\begin{theorem} \label{prop:SC}
There exists a  self-complementary $2$-complex without (tight, loose, or weak) Hamiltonian paths.
\end{theorem}

\begin{proof}
Consider the $2$-dimensional complex
    \[ S =  [1,3,5], [2,3,4], [1,2,4], [1,3,6], [1,2,5], [2,3,5], [1,4,5], [1,2,3], [1,2,6], [1,3,4]. \]
This $S$ is self-complementary: it is isomorphic to its complement 
            \[S^c = [1,4,6], [1,5,6], [2,3,6], [2,4,5], [2,4,6], [2,5,6], [3,4,5], [3,4,6], [3,5,6], [4,5,6]\]
via  the map $1\rightarrow 6$, $2\rightarrow 5$, $3\rightarrow 4$, $4\rightarrow 2$, $5\rightarrow 3$, $6\rightarrow 1$. 
Yet we verified with the software \cite{Pav21} that $S$ has no (tight, loose, or weak) Hamiltonian paths. 
\end{proof}

\begin{remark} As a curiosity, $S$ is also a counterexample to the generalization of other two well-known graph-theoretical statements: In fact, we verified with \cite{Pav21} that $S$ admits weakly Hamiltonian cycles but not weakly Hamiltonian paths, whereas all Hamiltonian graphs are traceable; moreover, $S$ is under-closed but not chordal, cf.~\cite{BSV}, whereas all interval graphs are chordal. 
\end{remark}

\begin{figure}[htb] 
    \centering
	\includegraphics[width=11em]{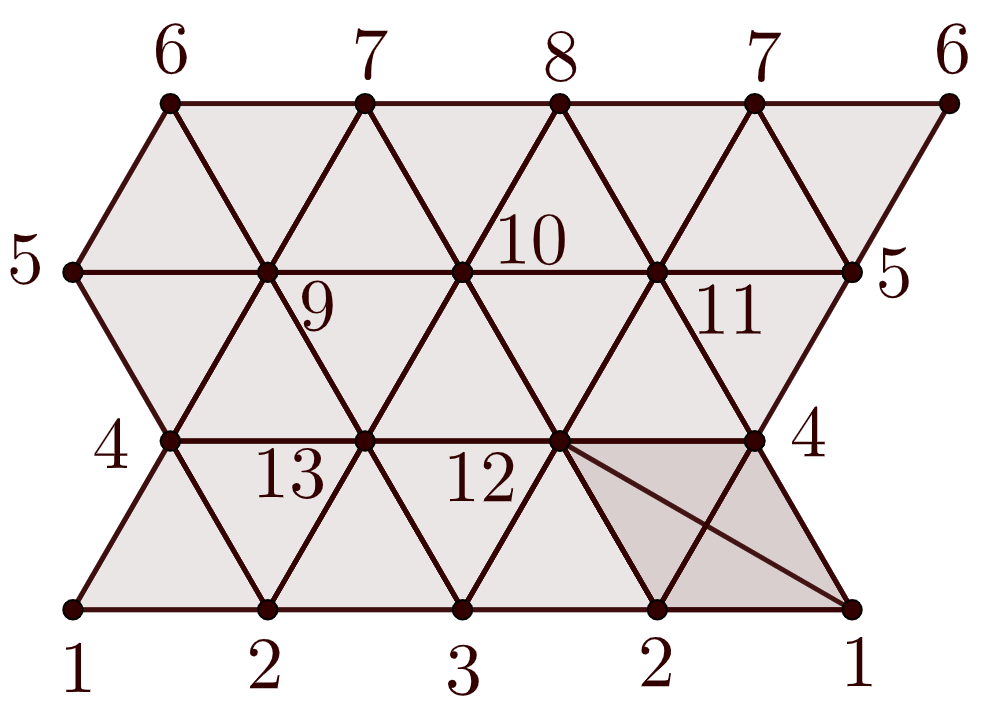}
	\caption{A 2-strongly-connected 2-complex $K$ whose square has no              Hamiltonian cycles.}
        \label{fig:GL12}
\end{figure}

Finally, recall that the \emph{square} $G^2$ of a graph $G$ is the graph obtained from $G$ by adding an edge $[i,j]$ for every two non-adjacent vertices $i,j$ of $G$ that are connected in $G$ by a $2$-edge path. Fleischner  showed that if $G$ is $2$-connected, then $G^2$ is Hamiltonian \cite{Flei}. 

\begin{definition} Given a $d$-complex $C$, let us define its \emph{square} $C^2$ as the $d$-complex obtained from $C$ by adding any $d$-face that would be introduced via any possible bistellar flip that does not add vertices. 
\end{definition}

For example, if $C=[1,2,3], [2,3,4], [3,4,5]$, then flipping at $[2,3]$ would introduce triangles $[1,3,4]$ and $[1,2,4]$, whereas flipping at $[3,4]$ would introduce $[2,3,5]$ and $[2,4,5]$; so $C^2$ is obtained by adding to $C$ the 4 triangles  $[1,2,4], [1,3,4], [2,3,5], [2,4,5]$.  Clearly, this notion of ``square of a $d$-complex'' boils down for $d=1$ to the square of a graph. 

\begin{theorem} \label{prop:Fl}
There exists a $2$-dimensional simplicial complex $K$ that is strongly-connected, and remains strongly-connected after the deletion of any vertex; yet its square $K^2$  does not admit (tight, loose, or weak) Hamiltonian cycles.  
\end{theorem}

\begin{proof} Consider the $2$-dimensional complex of Figure \ref{fig:GL12},
    \[ \begin{array}{ll}
      K =  &[2, 4,12], \; [1, 2, 4], \; [1, 2,12], \; [1, 4, 12], \; [2, 4, 13],\; [2, 3, 13], \;[3, 12, 13],\; [2, 3, 12], \\ 
    &[4, 5, 9], \; [4, 9, 13], \; [9, 10, 13], \;[10, 12, 13], \;[10, 11, 12], \;[4, 11, 12], \;[4, 5, 11],  \\
    &[5, 6, 9],  \; [6,7,9],\; [7, 8, 10], \;[7,8,11], \;[5, 6, 7], \;[5, 7, 11],\;[7, 9, 10], \;[8, 10, 11]. 
    \end{array}
     \]
    This $K$ is ``$2$-strongly-connected'', meaning that any complex obtained from $K$ by deleting at most one vertex, however chosen, is strongly-connected. The square $K^2$ is obtained by adding 46 more triangles. It can be verified with the software \cite{Pav21} that $K^2$ does not have (tight, loose, or weak) Hamiltonian cycles. 
\end{proof}



\end{document}